\documentclass[11pt]{article}

\usepackage{amsmath,epsf,amssymb,latexsym,cite}
 \usepackage{xypic}
 \usepackage{verbatim}
\usepackage{amsfonts}

\usepackage{fancyhdr}
\usepackage{indentfirst}
\usepackage{graphicx}
\usepackage{mathrsfs}
\usepackage{newlfont}
\usepackage[all]{xy}                                       

\usepackage{mathrsfs}
\usepackage{newlfont}
\usepackage{amssymb}
\usepackage{amsmath}
\usepackage{latexsym}
\usepackage{amsthm}
\theoremstyle{plain}                    
\newtheorem{teo}{Theorem}[section]      
\newtheorem{claim}[]{Claim}       
\newtheorem{lem}[teo]{Lemma}            
\theoremstyle{definition}               
\newtheorem{defin}{Definition}[section]
\theoremstyle{remark}                   

\renewcommand{\o}{\mathcal{O}}
\newcommand{\ox}{\mathcal{O}_X}
\renewcommand{\j} {\mathcal{J}}
\newcommand{\la}{\longrightarrow}

\newcommand{\q}{\mathbb{Q}}

\renewcommand{\k}{K_X}
\renewcommand{\j}{\mathcal{J}}

\begin{document} 
\title{Pluricanonical maps for threefolds of general type \thanks{The author would like to thank Professor Christopher Hacon for suggesting the problem and many useful conversations and suggestions.}}
\author{Gueorgui Tomov Todorov} 

\maketitle
\begin{abstract}In this paper we will prove that for a threefold of general type and large volume the second plurigenera is positive and the fifth canonical map is birational.
\end{abstract}
\section{Introduction.}  
 One of the main problems in algebraic geometry is to understand the 
structure of pluricanonical  maps.
If $X$ is a smooth complex projective variety  of general type and dimension $n$, then by definition 
the plurigenera $P_r=h^0(X,\ox(r\k))$ grow like $r^n$ and the 
pluricanonical maps $\phi_{r\k}:X\dashrightarrow \mathbb{P}(H^0(X,\ox(r\k)))$ are  birational 
for sufficiently large $r$. It is  natural to look for minimal 
values $r_n$ and $r_n'$ such that $P_{r}\neq 0$ for $r\ge r_n$ and  
$\phi_{r\k}$ is birational for $r \ge
 r_n'$. On a curve of general type the picture is clear, $P_r\neq 0$ for
 $r\ge 1$ and    $\phi_{r\k}$ is birational for $r \ge 3$. For surfaces we 
have that  $P_r\neq 0$ for $r\ge 2$ and    $\phi_{r\k}$ is birational for 
$r \ge 5$, by a  result of Bombieri  (cf. \cite{bom}).

 On threefolds the problem is much harder. The question has been extensively studied,  
for example  by Koll\'ar  \cite{kol3} and Benviste \cite{ben}. Chen and Hacon in \cite{chris3} 
have a proof that if the irregularity is positive then  $\phi_{7\k}$ is birational.
M. Chen \cite{chen} and S. Lee \cite{lee} prove that for $X$ minimal and smooth $\phi_{6\k}$ is birational.
By a recent result of J. Chen, M. Chen and D. Zhang in \cite{chenco} if  $X$ is a  minimal Gornenstein threefold   
 $\phi_{r\k}$ is birational 
for $r\ge 5$. However Tsuji \cite{tsu}, 
and later Hacon-M$^\textrm{c}$Kernan \cite{chris1}, Takayama \cite{tak} following ideas of Tsuji, have shown that  for a smooth 
projective variety of general type and dimension $n$ there exists an integer 
$r_n$, which depends only on $n$, such that   $\phi_{r\k}$ is birational 
for $r \ge r_n$. This theorem has also a nice consequence, namely if 
we fix a positive integer $M$ the family of smooth projective varieties
 of general type, dimension $n$ and volume smaller than $M$ is birationally bounded.

 Keeping this in mind it is natural to ask if we can find explicit values 
$r_n$ and $r_n'$ such that $P_{r}\neq 0$ for $r\ge r_n$ and  
$\phi_{r\k}$ is birational for $r \ge r_n'$ for a class of varieties of dimension $n$ and volume larger 
than a fixed positive number $M$. We will prove the following theorems:

 \begin{teo} Let $X$ be a projective threefold of general type   and assume that 
$$
\emph{vol} (X) > 2304^3.
$$

 Then $P_2=h^0(X,\ox(2\k))>0$.

 \end{teo}

\begin{teo} Let $X$ be a projective threefold of general type and assume that 
$$
\emph{vol} (X) > 4355^3.
$$

 Then $\phi_{5\k}$ is birational.

 \end{teo}
 Here the volume of  $X$  is defined as 
$$
\textrm{vol}(X)=\limsup_{m\to\infty}\frac{3!h^0(X,\ox(m\k))}{m^3}.
$$

 Note that the above result is optimal since  we can construct infinitely many
 families of threefolds of general type with increasing volume for which the $\phi_{r\k}$ is not birational 
for $r\le 4$. This can be done just by considering products of curves and surfaces of general type.

We will sketch the strategy of the proof for the non-vanishing of the second 
plurigenera. The proof of the second statement is similar. We would like to
 produce a divisor $D \sim \lambda \k$ such that $D$ has an isolated log canonical
 centre at a  general point $x\in X$. If $\lambda <1$, then we get a 
surjective homomorphism $H^0(X,\ox(2\k))\la H^0(X,\mathbb{C}_x)\cong\mathbb{C}$ and so $P_2 \ge 1$. 
Since the volume of $\k$ is large, by an easy dimension count, one finds a small rational number $\lambda>0$  
and a divisor   $D \sim \lambda \k$ which has a non-trivial log canonical centre at a general point $x\in X$. However
 it could happen that the dimension of the log  canonical center is positive. This 
is not a problem if the dimension of the log canonical centre at $x$ is one
 since then we can apply a standard technique to cut down the dimension of the log canonical centre 
 and have an isolated log canonical centre for 
$D' \sim \lambda' \k$ for some small positive rational number $\lambda'$ . If the dimension at $x$ of the log 
canonical centre is two we could still try to  cut down the dimension and indeed
 at the end we would produce a divisor  $D \sim \lambda \k$ with an isolated 
log canonical centre but $\lambda$ will not be small enough to give us a bicanonical
 section. So in this case a different approach is needed. In fact using ideas of 
M$^\textrm{c}$Kernan\cite{mac} and Tsuji we will produce a morphism to a curve and  we will use this morphism
to produce a section.

\section{Preliminaries.}        
\subsection{Notation.} 

We will work over the field of complex numbers $\mathbb{C}$. 
A $\q$-Cartier divisor $D$ is nef if $D\cdot C \ge0$ for any curve $C$ on $X$.
 We call two $\q$-divisors $D_1, D_2$ $\q$-linearly equivalent $D_1\sim D_2$ if 
there exists an integer $m>0$ such that $mD_i$ are Cartier and linearly equivalent.
We call two $\q$-Cartier divisors $D_1, D_2$ numerically equivalent $D_1\equiv D_2$ if $(D_1-D_2)\cdot C=0$  for any curve $C$ on $X$.
A log pair $(X,\Delta)$ is a normal variety $X$ and an effective $\q$-Weil 
divisor $\Delta$ such that $\k+\Delta$ is $\q$-Cartier. A projective morphism 
$\mu:Y \la X$ is a log resolution of the pair  $(X,\Delta)$ if $Y$ is smooth and 
$\mu^{-1}(\Delta)\cup\{\textrm{exceptional set of } \mu\}$ is a divisor 
with simple normal crossing support. For such $\mu$ we write 
$\mu^*(\k+\Delta)  =K_Y+\Gamma$, and $\Gamma=\Sigma a_i\Gamma_i$ 
where $\Gamma_i$ are distinct integral divisors. A pair is called 
klt (resp. lc) if there is a log resolution $\mu:Y \la X$ such that in 
the above notation we have $a_i <1$ (resp. $a_i\le 1$). The number $1-a_i$ is called 
log discrepancy of $\Gamma_i$ with respect to the pair  $(X,\Delta)$.  
We say that a subvariety $V \subset X$ is a log canonical centre if it is the 
image of a divisor of log discrepancy at most zero.
 A log canonical place is a valuation corresponding to a divisor of log discrepancy at most zero.
 A log canonical centre is pure if $\k+\Delta$ is log canonical at the generic point of $V$. If moreover there is a unique log canonical place lying over the generic point of V, then we say that $V$ is exceptional. 
LCS$(X,\Delta,x)$ is the union of all log canonical centres of $(X,\Delta)$ through the point $x$.
We will denote by LLC$(X,\Delta,x)$ the set of all log canonical centres 
containing a point $x \in X$. 

\subsection{Volumes.}
\begin{defin} Let $X$ be irreducible projective variety of dimension $n$ and let $D$ be a 
$\q$-Cartier divisor. The volume of $D$ is 
$$
\textrm{vol}(D)=\limsup_{m\to\infty}\frac{n!h^0(X,\ox(mD))}{m^n}.
$$
\end{defin} 

 When $D$ is nef we have that vol$(D)=D^n$, so that for a very ample divisor $D$ the volume of $D$ is just the degree of 
the image of $X$ given by the linear system $|D|$ (cf. \cite{raz}, Section 2.2C). The volume depends only on 
the numerical class of $D$ and it can be extended to a continuous function 
vol$:N^1_{\mathbb{R}}\la{\mathbb{R}}$, where
 $N^1_{\mathbb{R}}=N^1(X)\otimes\mathbb{R}$ and $N^1(X)$
 is the Neron-Severi group of $X$(cf. \cite{raz}, Theorem 2.2.43).
 The volume is invariant under birational transformations(cf. \cite{raz}, Proposition 2.2.42).

 The following lemma is a standard tool to produce log
 canonical centres at a general point of a variety.

\begin{lem}\label{mult}(cf. \cite{raz}, Lemma 10.4.11) Let $X$ be an  irreducible projective variety of dimension $d$,
 $L$ a $\q$-Cartier divisor and $x \in X$ a smooth point. If for some positive rational number $\alpha$ we have that
$$
\emph{vol} (L) > \alpha^d,
$$
then for any sufficiently divisible integer $k\gg 0$ there exists a divisor 
$$
A=A_x \in |kL|\qquad \textrm{with} \qquad \emph{mult}_x(A)>k\alpha.
$$
Moreover we can take $k$ to be independent of the smooth point.
\end{lem}

 We will also use the following variant of Lemma \ref{mult}.

\begin{lem}\label{mult1} Let $X$ be an irreducible projective variety of dimension $d$,
 $L$ a $\q$-Cartier divisor and $x, y \in X$ smooth points. 
If for some positive rational number $\beta$ we have that
$$
\emph{vol} (L) >2 \beta^d,
$$
then for any sufficiently divisible integer $k\gg 0$, there exists a divisor 
$$
A=A_{x,y} \in |kL|\qquad \textrm{with} \qquad \emph{mult}_x(A)>k\beta\textrm{ and \emph{mult}}_y(A)>k\beta.
$$
\end{lem}

\begin{proof} The prove of this fact is well known. For a similar argument see for example \cite{raz}, Proposition 1.1.31.
\end{proof}
\subsection{Multiplier Ideals.}

\begin{defin} Let  $(X,D)$ be a log pair with $X$ smooth  and  let $\mu : Y \la X$ be a log resolution. Then define the multiplier ideal sheaf of $D$ to be 
$$
\j(X,D)=\mu_*\o_{Y}(K_{Y/X}-\llcorner\mu^*D\lrcorner)\subset\ox.
$$
When no confusion is  possible we will write simply $\j(D)$.
\end{defin}
 We will recall some facts, which will be needed in what follows.  We will repeatedly make use of the following 
generalization of the Kodaira vanishing theorem, which will be referred to as Nadel vanishing. For a proof of this theorem we refer to  section 9.4.B of \cite{raz}.
\begin{teo}\label{nadel} Let $X$ be a smooth projective variety, 
$D$ a $\q$-divisor  and $L$ be an integral divisor such that $L-D$ is nef and big.
Then 
$$
H^i(X,\ox(\k+L)\otimes\j(D))=0\qquad \textrm{for} \qquad i>0.
$$
\end{teo}
 
Define the log canonical 
threshold of $D$ at $x$ to be  
$$
c(D,x)=\inf\{ c>0|\j(c\cdot D)\textrm{ is non-trivial at }x\}.
$$

 We will also define 
$$
\textrm{Nklt}(X,D)=\textrm{Supp}(\ox/\j(X,D))\subset X
$$
with the reduced structure and we will call it the \emph{non-klt locus} for $(X,D)$ 
 
Also we will say that $(X,D)$ is \emph{klt} (resp. \emph{lc}) at $x\in X$ 
if $(U,D_{|U})$ is klt (resp. lc) for some Zariski open neighborhood $U$ of $x$.
 It is often useful to assume that LCS$(X,D,x)$  
is irreducible at $x$. The next lemma asserts that one can achieve this after an  arbitrarily 
small perturbation of $D$.

\begin{lem}\label{irr}(cf. \cite{chris1} Lemma 2.5, \cite{kaw1}, \cite{ambro} ) Let $X$ be a smooth projective variety and $\Delta$ an effective $\q$-divisor and assume that $(X,\Delta)$ is lc at some point $x\in X$. If $W_1,W_2\in$LLC$(X,\Delta,x)$ and $W$ is an irreducible component of $W_1\cap W_2$ containing $x$, then $W\in$LLC$(X,\Delta,x)$. Therefore, if $(X,\Delta)$ is not klt at $x$, then LLC$(X,\Delta,x)$ has a unique minimal irreducible element, say $V$. Moreover, there exists an effective $\q$ divisor $E$ such that 
$$
\textrm{LLC}(X,(1-\epsilon)\Delta+\epsilon E,x)=\{V\}
$$
for all $0<\epsilon \ll 1$. We may also assume that there is a unique place laying above $V$ and if $x\in X$ is general and $L$ is a big divisor, then one can take $E=aL$ for some positive number $a$.
\end{lem}

 The next theorem shows how multiplier ideal sheaves behave with respect to restrictions to  smooth divisors.
\begin{teo}\label{inclu}(cf. \cite{raz}, Theorem 9.5.1)   Let $X$ be a smooth projective variety, let $D$ be an effective 
$\q$-divisor on $X$ and let $H\subset X$ be a smooth irreducible hypersurface which is
 not contained in the support of $D$. Then there is an inclusion
$$
\j(H,D_{H})\subset \j(X,D)_H.
$$
\end{teo}
Here $D_H$ is the restriction of $D$ to $H$, which is by assumption an effective $\q$-divisor on $H$  so that 
the multiplier ideal sheaf in question is defined. $\j(X,D)_H$ indicates the image of $\j(X,D)$ under the natural maps
$$
\j(X,D)\hookrightarrow \ox\la\o_H.
$$

Finally we will use the following Lemma from Hacon-M$^\textrm{c}$Kernan (\cite{chris1}, Lemma 2.6).
\begin{lem}\label{sec} Let $X$ be a smooth projective variety and $D$ a big divisor on 
$X$. Let $x$ be a general point, $\lambda$ a rational number with  $0<\lambda <1$ and $\Delta\sim\lambda D$ a $\q$-divisor such that $x$ is a component of $\emph{LCS}(X,\Delta,x)$. Then $h^0(X,\ox(\k+D) )>0$.
\end{lem}
\section{Proof of the Theorem 1.1.}

 Let us start by replacing $X$ by a desingularization and  assume   that 
$$
\textrm{vol}(X) > \alpha ^3.
$$
The value of $\alpha$ as in Theorem 1.1 will be determined in the course of the proof.
 We are going to proceed in the following way. The first step 
will be to produce a log canonical centre at some point $x\in X$
 and then to cut  down the dimension of the log canonical centre 
until it is an isolated point. This can be done as long as 
the log canonical centre is not a surface of small volume. In this
 last case we are able to produce a morphism to a curve with 
special properties which will help us to produce sections of multiples of 
the canonical bundle.

 Our main tool to cut down the dimension of the log canonical centres 
is the following result of Hacon-M$^\textrm{c}$Kernan \cite{chris1}. Fix a log pair
$(X,\Delta)$ where $X$ is $\q$-factorial and $\Delta$ is an effective 
$\q$-divisor. Let $V$ be an exceptional log canonical centre of $\k+\Delta$. 
Let $f:W\la V$ be a resolution of $V$ and let $\Theta$ be a $\q$-divisor on $W$.
 Suppose that there are positive rational numbers $\lambda$ and $\mu$ such 
that $\Delta\sim \lambda\k$ and $\Theta \sim \mu K_W$. Let $\nu=(\lambda+1)(\mu+1)-1$. 
We suppose that $V$ is not contained in the augmented  base locus of $K_X$ and $W$ is 
of general type. Here the augmented  base locus is defined as follows. For a divisor $L$ by the stable base locus {\bf B}$(L)$ we mean the algebraic set that is 
 $\cap\textrm{Bs}(|mL|)$, where the intersection is taken over all positive integers $m$. For a big divisor $L$ the sets {\bf B}$(L-\epsilon A)$ are the same for any ample divisor $A$ and any $0<\epsilon\ll 1$ (cf. \cite{raz}, Lemma 10.3.1) and this is the augmented base locus of $L$.
  In this setting we have the following result.

\begin{teo}\label{cut1}(Hacon-M$^\textrm{c}$Kernan) There is a very general subset $U$ of 
 $V$ with the following property. Suppose that $W'\subset W$ is a pure log  
canonical centre of $K_W + \Theta$, whose image $V' \subset V$ intersects U. 
Then for every positive rational number $\delta$ we may find a divisor $\Delta'$
 on $X$ such that $V'$ is an exceptional log canonical centre of $\k+\Delta'$ where
 $\Delta' \sim (\nu+\delta)\k$.

\end{teo}

\indent \emph{Proof of Theorem 1.1}.
Denote by  $X_0$ the complement of all subvarieties contained in $X$ that are not of general type 
(since $X$ is of general type, this is a countable union of closed subsets) and of the augmented base locus of $\k$. 
Choose a general point $x \in X_0$.

 Since vol$(\k)>\alpha^3$  by Lemma \ref{mult}   
for any  integer $k \gg 0$ there exists a divisor 
$$
A = A_x \in |k \k| \qquad \textrm{with} \quad \textrm{mult}_x(A)> k\alpha.
$$

So, if we consider the divisor 
$$
\Delta ' = \lambda ' \frac{A}{k}
$$
for $\lambda' < \frac{3}{\alpha}$ but close enough to $ \frac{3}{\alpha}$, we have that 
$$
\textrm{mult}_x \Delta ' >3.
$$ 
 But then $(X,\Delta')$ is not klt at $x$ and so by Lemma \ref{irr}, if $V$ is the unique minimal 
irreducible element of LCC$(X,\Delta',x)$, we have that 
$$
\textrm{LLC}(X,(1-\epsilon)\Delta'+\epsilon E,x)=\{V\}
$$
for all $0<\epsilon\ll0$. Here we can take $E\sim a\k$ for some positive integer $a$. 
Now by taking  $\epsilon < \frac{3\alpha^{-1}-\lambda'}{a-\lambda'}$ we obtain a 
divisor $D$ such that $D\sim \lambda \k$ with $\lambda<\frac{3}{\alpha}$ and 
$$
\textrm{LCC}(X,D,x)=\{V\}.
$$
We will analyze different cases based on the dimension of this log canonical centre.

\emph{Case I - \emph{dim}$V$=0.}

To produce a section of $H^0(X,\ox(2\k))$ we want to apply 
Lemma \ref{sec}. So what we need is that  $\lambda < 1 $. 
This can be achieved  for $\alpha =3 $.

 It is now clear that $P_2(X) >0$.

\emph{Case II - \emph{dim}$V=1$.}

 We want to cut down the dimension of the log  canonical centre  to reduce 
to the previous case. We wish  to apply Theorem \ref{cut1}, 
so we need to produce a divisor of high multiplicity at a very general smooth point of $V$.
 Let $f:W\la V$ be the normalization of $V$. By our choice of $x$ we know that $W$ is of general type, so $\textrm{vol}(K_W)\ge 2> \frac{15}{8}$. 

 Lemma \ref{mult} implies that for a smooth point $w \in W$ there exists a $\q$-divisor $\Theta' \sim K_W$ with 
multiplicity greater than $\frac{15}{8}$ at $w$. Then $\Theta=\mu\Theta'\sim\mu K_W$ has multiplicity one  at $w$ for some $\mu<\frac{8}{15}$.  By Theorem \ref{cut1}
 there is a divisor $D'$ on $X$ such that
$$D' \sim (\lambda+\mu+\lambda\mu+\delta)\k,$$ 
where $\delta$ is any small positive rational number, and the pair  $(X,D')$ has some very general point of $V$ as  an exceptional log canonical centre of singularities.
 Choose $\delta=10^{-6}$.
 As before we just have  to make sure that $\lambda+\mu+\lambda\mu+\delta <1$. Now 
$\lambda+\mu+\lambda\mu+\delta< (\frac{8}{15}+\frac{23}{15}\lambda+\delta)$ so it is enough
 that  $\frac{8}{15}+\frac{23}{15}\lambda+\delta<1$,
and since $\lambda<\frac{3}{\alpha}$ this is possible  
 as long as 
$\alpha \ge 10$. 

\emph{Case III - \emph{dim}$V=2$}.

Let $f:W \la V$ be a resolution of singularities of $V$. As in the previous case $W$ is of general type.
We will divide our argument in two cases based on whether  the volume of $K_W$ is large or small.
In the first case when the volume of $K_W$ is large we may  proceed as before, that is  produce a log canonical centre of 
$t\k$ for some $t<1$ and of smaller dimension. In the small volume case  we will show that there exists a fibration 
onto a curve and using this fibration we will be able to prove the existence of the required sections.  

Suppose now that  $\textrm{vol}(K_W) > 64$. 
By Theorem \ref{cut1} and   using that $\textrm{vol}(K_W) > 64$ we can guarantee there is a divisor $D'\sim (\frac{1}{4}+\frac{5}{4}\lambda+\delta)\k$ such that $V'$ is an exceptional log canonical centre of $(X,D')$ of dimension 
  not greater than one. 
If $\textrm{dim} V'=0$ then we are done by the discussion in Case I. If instead $\textrm{dim} V'=1$ replacing $V$ by $V'$ and using the discussion in Case II we see
that  there is a divisor $D''$ such that   
$$
D''\sim \left(\frac{8}{15}+\frac{23}{15}\left(\frac{1}{4}+\frac{4}{5}\lambda+\delta\right)+\delta\right)\k\sim \left(\frac{55}{60}+\frac{92}{75}\lambda+\frac{38}{15}\delta\right)\k
$$
for which for LLC$(X,D'',x')=\{x'\}$ for a point $x'\in X$. As we have seen already, once we have an isolated centre  we are done as long as $\frac{55}{60}+\frac{92}{75}\lambda+\frac{38}{15}\delta<1$, and this can be achieved if $\alpha > 45$.
 
 To proceed we need a simplified version of a Lemma 3.2 from M$^\textrm{c}$Kernan \cite{mac}. 
For the convenience of the reader we include a proof  of the result. Recall that a subset $P$ of 
$X$ is called countably dense if it is not contained in the union of countably many closed subsets of $X$.
\begin{lem}\label{mac} Let $Y$ be a smooth projective variety and suppose that for every point $y\in P$, where $P$ is a countably dense subset of $Y$, we may
find a pair $(\Delta_y,W_y)$ such that $W_y$ is a pure log canonical centre for $K_Y+\Delta_y$ at $y$ and 
$\Delta_y \sim \Delta /w_y$ for some big Cartier divisor $\Delta$ on $Y$. Then there exists a diagram 
$$
  \xymatrix{
   {Y'}\ar[r]^{\pi} \ar[d]^{f}
   & Y\\
   B&\\
   }
$$
such that $f$ is a dominant morphism of  normal projective varieties with connected fibres and for a general fibre $Y'_b$ of $f$ 
there exists $y \in f(Y'_b)$ so that $f(Y'_b)$ is a pure log canonical centre for $K_Y+\Delta_y$
with $\Delta_y\sim\Delta/w$ at $y$, for  some $w$.
Also   $\pi$ is a generically finite and dominant morphism of normal varieties.
\end{lem} 

\begin{proof}
As $\Delta$ is big by Lemma \ref{irr}, possibly passing to a subset of $P$, 
we may assume that $W_y$ is an exceptional log canonical centre. If we decompose a countably
 dense set as a countable union at least one of the sets in the union is
countably dense. The field of rational numbers is countable so on a 
countably dense subset $w_y$ does not depend on $y$. 
 We obtain a countably dense subset $Q$ of $P$ such that $w_y=w$, for some fixed rational number $w$,
for every $y\in Q$.
  Then there is an integer $m'$ such that $m'\Delta_y$ and $m'\Delta /w$ 
are integral and linearly equivalent, so $m'\Delta_y\in |m'\Delta /w|$. 
 In the linear system  $|m'\Delta/w|$  take 
the closure $B$ of the points that correspond to divisors $\Delta_y$.
 Replace $B$ by an irreducible component which contains a set corresponding to every 
point of a countably dense subset of $Y$ and let $f:H\la B$ be 
the universal family over $B$, $H\subset B \times Y$. Pick a log resolution of the generic 
fiber and extend this to an embedded resolution of $Y \times U$ for some open subset $U$ of $B$.
 By our assumptions there 
is a unique exceptional divisor of log discrepancy zero over the generic 
point of the base corresponding to the log canonical centres $W_y$. 
Thus by taking a finite cover of $B$ and 
passing to an open subset of $U$ we may assume that there is a
 morphism $f:Y'\la B$ whose fibre over $b\in U$ is a log canonical centre for $K_Y+\Delta_b$. 
 
 By passing to an open subset  of $U$  we may assume that $f$ is flat, 
so $U$ maps to the Hilbert scheme. Replace $B$ by the normalization 
of the closure of the image of $U$ in the corresponding  Hilbert scheme and $Y'$ 
by the normalization of the pull back of the universal family. 
Cutting $B$ by hyperplanes we may assume that the map from $Y'$ to $Y$  is generically finite.
\end{proof}

 Going back to the proof of Theorem 1.1 we have left to discuss the case in which for every very general
$x \in X$ we have a pair of $(D_x,V_x)$, such that
\begin{description}
\item[(1)] $D_x \sim \lambda \k$.
\item[(2)] $V_x$ is a pure log canonical center of $D_x$.
\item[(3)] $\dim V_x=2$.
\end{description}
Observe that here we can take the same $\lambda$ for every point in a countably dense subset of $X$ by the argument of the first part of the proof of the Lemma above. Now by the  previous lemma we have a diagram:
 $$
  \xymatrix{
   {X'}\ar[r]^{\pi} \ar[d]^{f}
   & X\\
   B&\\
   }
$$
where $\pi$ is dominant and generically finite of normal varieties,
 and the image of the general fiber of $f$ is  $V_x$ for some $x\in X$. 

 We claim that we can also assume that the map $\pi$ is birational. There is a general proof of this fact in  M$^\textrm{c}$Kernan's paper \cite{mac}.
In our case, because of the low dimension, we present a simpler argument. We will use an elementary lemma of calculation of centres 
of log canonical singularities.
\begin{lem}\label{centre} Let $Y$ be a smooth projective variety and let $(\Delta_i, W_i)$ for $i=1,2$ be a pair such that at some 
point $y\in Y$ we have that $W_i$ is an exceptional centre of log canonical singularities of codimension $1$ at $y$ 
for $K_Y+\Delta_i$ with $\Delta_i$ smooth at $y$. Then there exists $Z\subset W_1 \cap W_2$  a minimal pure centre of log canonical singularities at $y$  for a pair $(Y,\Delta)$  with 
$\Delta=k(\Delta_1+\Delta_2)$ for some rational number $0<k\le1$.
\end{lem} 
 \begin{proof}

 First of all notice that if we set $c:=c(Y,\Delta_1+\Delta_2,y)$  we have that $c\le 1$. Also LCS$(X,c(\Delta_1+\Delta_2),y)\subset W_1 \cap W_2$.
The lemma now follows from the first part of Lemma \ref{irr}.
\end{proof}

 
 Let $X_b$ be the fibre of $f$ over $b\in B$. We claim that there exists a closed subset $X_1$ of $X$ such that  for any pair $(D_x,V_x)$ as above, $x\notin X_1$ and $\pi(X_b)=V_x$ for some $b\in B$ we have that $D_x$ is  smooth at $x$. In fact there is an open subset $U$ of $B$ that is parametrizing a universal family $Y$ of divisors $D_b$ so that for $x\in \pi(f^{-1}(U))$ and $f(X'_b)=V_x$ we have $D_x=D_b$. But after shrinking $U$ the singular locus of $Y$, which is at most two dimensional, dominates the union of singularities of divisors $D_x$ with  $x\in \pi(f^{-1}(U))$.  Since $B\backslash U$ is a finite set the claim follows. 


 Suppose now that   the inverse image of a general $x\in X\smallsetminus X_1$ under $\pi$ is contained in two different fibres of $f$. Then  through some general point $x\in X$ there are two log canonical centres, say  $V_1$ and $V_2$, of the pairs  $(X,D_1)$ and $(X,D_2)$,
with $D_i \sim \lambda \k$ for $i=1,2$.
 But since $x\notin X_1$ it follows that $D_i$ are smooth at $x$ for  $i=1,2$ and 
we can apply the previous lemma. So there is   $C$  a minimal  pure log canonical
 centre of $(X,D)$ with $D=k(D_1+D_2)$ and $0<k\le 1$ contained in  the intersection of $V_1$ and $V_2$ through
 $x$. 
We have  $D\sim  k(2\lambda) \k$ . Since dim$C\le 1$ by the analysis in Case I and II 
if $\frac{8}{15}+\frac{23}{15}(2\lambda)<1$ we obtain a section of $2\k$. The inequality can be arranged for $\alpha >20$.

 So we will replace our $X$ with an appropriate smooth birational model such  that we have a morphism to a curve $f: X \la B$ and  for a general fiber 
$X_{b'}$ over a point $b'\in B$ there is a divisor $D_{b'}\sim \lambda \k$ for which we have 
$$ 
\j(D_{b'})\subset\ox(-X_{b'}).
$$  

 Now run the relative minimal model program for the morphism $f:X\la B$. 
Denote the resulting morphism by $f_0:X_0\la B$ and here $X_0$ is $\q$-factorial with terminal singularities. 
Since the singularities of $X_0$ are isolated points (cf. \cite{kol4}), we can replace $X_0$ by a resolution which does not change the general fibre of $f_0$. Thus by replacing $f$ by the appropriate resolution we can assume that the general fibre of $f$ is a minimal and  smooth surface. 
  
We need the following lemma.

 \begin{lem}\label{vanish}  Let $g:Y\la Z$ be a surjective morphism from a smooth projective variety  $Y$ 
to a normal variety $Z$. Let $L$ be a divisor on $Y$ such that 
$L \equiv g^*M + \Delta$, where $M$ is nef and big $\q$-Cartier divisor on $Z$ and $\Delta$ is an effective $\q$-Cartier divisor on $Y$.
 Then $H^i(Z,R^jg_*(\o_Y(K_Y+L) \otimes \j(\Delta))=0$ for $i >0$, $j\ge 0$.

\end{lem}
\begin{proof} Let $g':Y'\la Y$ be a log resolution of the pair $(Y,\Delta)$ and let $h=g\circ g'$. 
By Corollary 10.15 of \cite{kol1} we have that 
$H^i(Z,R^jh_*(K_{Y'}-\llcorner g'^*\Delta\lrcorner+g'^*L))=0$ for $i >0$, $j\ge 0$. Notice that we have 
$R^jg'_*\o_{Y'}(K_{Y'/Y}-\llcorner g'^*\Delta\lrcorner)=0$ for $j>0$ (cf. \cite{raz}, Theorem 9.4.1) and so the Lemma follows by the projection formula. 
\end{proof}

For an integer $m$ sufficiently large and divisible consider a general divisor $G\in |m(\k+X_b)|$ for some general $b\in B$.
For a general $b'\in B$ set $D=D_b'$ and let $F=3D+\frac{1-3\lambda}{m}G+(3\lambda-2)X_{b'}$. We have that $\j(F)\subset\ox(-X_{b'})$. 
The divisor $F$ is effective since $X_{b'}$ is a pure log canonical centre for $(X,D)$.
Consider 
$$
E:=f_*(\o_{X}(2K_{X})\otimes\j(F)).
$$
This is a vector bundle, since it is  a torsion free sheaf on a curve.


\begin{claim} Let  $r$ be the rank of the bundle $E$. If $\alpha > 2304$ then we have that $r>0$.
\end{claim}
\indent \emph{Proof of the Claim.}
   The  general fiber $X_b$ is
 a smooth surface and it  is enough to 
prove that 
$$
H^0(X_b,\o_{X_b}(2K_{X_b})\otimes\j(F)_{|X_b})\neq 0.
$$
 The hypothesis of the restriction Theorem \ref{inclu} for multiplier ideal sheaves are satisfied, so we have that 
 $\j(X_b,(3D+\frac{1-3\lambda}{m}G)_{|X_b})\subset\j(X,3D+\frac{1-3\lambda}{m}G)_{|X_b}$. 
Also since $\k$ is big we 
can apply Kawamata's Theorem A from \cite{kaw3}, so we obtain 
a surjection $H^0(X,\ox(m(K_X+X_b)))\twoheadrightarrow 
H^0(X_{b},\ox(mK_{X_b}))$. But  $X_b$ is minimal, so  $|mK_{X_b}|$ is base point free. This implies that 
$\j(X_b,(3D+\frac{1-3\lambda}{m}G)_{|X_b})=\j(X_b,3D_{|X_b})$. Define  $\Delta=D_{|X_b}$ and observe that since 
$(K_X)_{|X_b}\sim K_{X_b}$ we have that $\Delta\sim \lambda K_{X_b}$. So it is enough to show that 
$$
H^0(X_b,\o_{X_b}(2K_{X_b})\otimes\j(3\Delta))\neq 0.
$$

 By a result of Bombieri (cf. \cite{bom}) we have that the linear system $|4K_{X_b}|$ is base point free. 
If we arrange that $3\Delta\cdot C<1$, for a curve $C$ that is 
a component of a divisor in  $|4K_{X_b}|$ and not a component of $\Delta$, we will have that mult$_x3\Delta<1$ for $x\in \Delta$, and this  implies that $\j(3\Delta)$ is trivial (cf. \cite{raz}, Proposition 9.5.11). The inequality could be arranged for $\lambda < \frac{1}{768}$ and thus for $\alpha>2304$.  But the second plurigenera on a surface of general 
type is always non-zero, and the claim follows.
\qed

 By Lemma \ref{vanish} since $\k\equiv X_b+F$ we have the vanishing  
$$
H^1(B,f_*(\ox(2 \k)\otimes\j(F))=0.
$$  
Tensoring  the exact sequence 
$$
0\la\o_B\la\o_B(b')\la\o_{b'}\la 0
$$ 
by $E$ we get a section of $f_*(\ox(2 \k)\otimes\j(F))\otimes\o_B(b')$. But  we have the inclusion
$$f_*(\ox(2\k)\otimes\j(F))\otimes\o_B(b') \hookrightarrow f_*(\ox(2\k))$$
and  we obtain  a section of $\ox(2\k)$.

\qed

\section{Proof of Theorem 1.2.}

We want to prove that the map given by the linear system $|5\k|$ is birational. 
So we would like to produce a divisor $D \sim t\k$  with $x,y \in$LLC$(X,D)$, where $t$ is a positive rational number. Then  apply Kawamata-Viehweg vanishing to lift sections from
 $\mathbb{C}_x\oplus\mathbb{C}_y$. Unfortunately if we proceed as above
 we are not able to do so for $t<4$. We will have to use a slightly different approach.

 As before we will divide the proof in two parts. 
First we will study the case when all the surfaces that appear as
 log canonical centres of some divisor linearly equivalent to a multiple
 of the canonical class  are of large volume. In this case we deduce the result from   Takayama's results  \cite{tak} by carefully choose the necessary constants.

\indent \emph{Proof of Theorem 1.2}. Start by replacing $X$ by a disingularization.  Recall that we are assuming 
$$
\textrm{vol}(\k)>\beta^3.
$$

By Theorem 3.1 and  Proposition 5.3 (see also Notation 5.2) of Takayama \cite{tak} for  $\epsilon>0$ there is a birational morphism $\mu :X'\la X$  and a decomposition
$$
\mu^*\k\sim A+M,
$$
where $A$ is an ample $\q$-divisor and $M$ is an effective divisor whose support contains all the exceptional locus of $\mu$, such that  for two general points $x_1$ and $x_2$ there is  a divisor $D$,
$$
D\sim \lambda A
$$
with 
$$
\lambda<
\lambda'_1+\frac{4\sqrt{2}}{(1-\epsilon)\alpha_1}+2\frac{2+\lambda'_1+\frac{4\sqrt{2}}{(1-\epsilon)\alpha_1}+2\epsilon}{(1-\epsilon)\alpha_2}+2\epsilon, 
$$
such that $x_i\in$Nklt$(X',D)$ for $i=1, 2$ and codim Nklt$(X',D)=3$ around at least one of the $x_i$, $i=1,2$. 
 Then by  Nadel vanishing we see that $(\ulcorner \lambda \urcorner+1)\k$ gives a birational map. Here
$\alpha_1$ is the square root of  a lower bound for the volume of the surfaces of general type contained in $X$,
$\alpha_2$ is a lower bound for the volume of the curves of general type contained in $X$, and 
$$
\lambda'_1<2^{1/3}\frac{3}{(1-\epsilon)\beta}.
$$

In any case we would like to ensure that  $\lambda <4$. 
 First note that for any curve $C$ of general type we have $\textrm{vol}(C)\ge 2$. So in the inequality above we
may  take 
$\alpha_2=2$. We also may  choose $\epsilon=0.01$. 
Since we will argue as before we will take $\alpha_1=6$. To achieve that  $\lambda <4 $ it suffices to be able to choose $\lambda_1' <\frac{1}{20}$ , and this can be done as long as   $\beta > 76$.

 Finally arguing as in the proof of Theorem 1.1 and the proof of Lemma \ref{mac} we can assume that through every point $x\in X$ in a
countably dense subset of $X$ we have a surface $V_x$ 
that is a centre  of log canonical singularities for  $(X,D_x)$, $D_x \sim \lambda_1' \k$ and so we get a diagram 
 $$
  \xymatrix{
   {X''}\ar[r]^{\pi} \ar[d]^{f}
   & X\\
   B&\\
   }
$$
 As before we argue that we can assume that $\pi$ is birational. In fact let $x$ and
 $y$ be two general points of $X$ as above. We would like to produce a divisor $D$ that is lc, 
but not klt at one of the points with centre a curve and that is not klt at the other point and apply the 
inductive steps of Takayama. Consider the divisor $D_x+D_y$. By rescaling we can assume that 
we have a divisor that is lc, but not klt at one of the points, say $x$, and it is not klt at $y$. That is, we
 have $D_1\sim \lambda' \k$ with $\lambda'\le2\lambda_1'$. Also we may assume that the
 centre of log canonical singularities at $x$ is a surface. If $\pi$ is not birational then 
there is another divisor $D_x'$ with centre at $x$ another  surface. Consider $D'=D_1+kD_x'$,
 where $k=\max\{c,\textrm{such that }(X,D_1+cD_x')\textrm{ is lc at }x\}$. Since $x$ and $y$ are general points 
we observe that $0<k\le1$.
  By possibly applying Lemma \ref{centre},
 we obtain a divisor $D\sim \lambda''\k$ with the desired properties and   $\lambda''\le 3\lambda_1'$.
Then Proposition 5.3 of \cite{tak} implies that by taking 
$$
\lambda<3\lambda_1'+\frac{2+3\lambda_1'+\epsilon}{1-\epsilon}+\epsilon
$$ 
we have that $(\ulcorner \lambda \urcorner+1)\k$ gives a birational map. We can take
$\lambda<4$ for $\beta\ge 16$. 

Thus we may assume that we have a morphism to a curve $f: X \la B$, such that the  general fiber 
$X_b$ of $f$ is minimal and smooth and it is log canonical centre for $(X,D_x)$
 with $D_x \sim \lambda \k$ and $\lambda<2^\frac{1}{3}\frac{3}{(1-\epsilon)\beta}$.
Moreover we have that the vol$(X_b)<36$ .

%

 The goal is to find a section of $\ox(5\k)$ that separates two general points. In the case when the  two general points  $x$ and $y$ in $X$ do not lie in  the same fiber, say $x\in X_{b_1}$ and $y\in X_{b_2}$, 
we will create a section of $f_*\ox(5\k)$ that separates $f(x)$ and $f(y)$. For $b$ and $b'$ general points in $B$ 
consider a general divisor  $G\in|m(\k+X_{b'})|$ for some integer $m$
sufficiently large and divisible and define $F=D_{b_1}+D_{b_2}+6D_b+\frac{4-8\lambda}{m}G+(8\lambda-5)X_b$.
 Then  
arguing as in Claim 1 we can guarantee that 
 $E=f_*(\ox(5\k)\otimes\j(F))$ is a vector bundle of positive rank for $\lambda < \frac{1}{1152}$ so that $\beta 
\ge 4355$.  Since $4\k\equiv X_b+F \sim f^*b+F$ we can apply Lemma \ref{vanish}. It follows that  $H^1(B,f_*(\ox(5\k)\otimes \j(F)))=0.$ Then tensoring  the short exact sequence 
$$
0\la\o_B\la\o_B(b_1+b_2)\la \o_{b_1}\oplus\o_{b_2}\la0,
$$
with $E$  we see that  there is  a section of  $f_*(\ox(5\k)\otimes\j(F))\otimes\o_B(b_1+b_2)$ that separates $b_1$ and $b_2$.  But then the inclusion
$$
 f_*(\ox(5\k)\otimes\j(F))\otimes\o_B(b_1+b_2)\hookrightarrow f_*\ox(5\k)
$$
implies that the linear system $|5\k|$ separates points that are not on the same general  fibre of $f$.
 
 Suppose now that $x$ and $y$   lie in  the same fiber, say $X_b$ over $b\in B$. Let $b'$  be another general point and let $(D_{b'},X_{b'})$ be the corresponding pair such that  $\j(D_{b'})\subset\ox(-X_{b'})$.  Write $\k\sim A'+M'$ as a sum of an ample and an effective divisor. Consider a general divisor  $G\in|m(\k+X_b)|$ for some integer $m$
sufficiently large and divisible and define $F=6D_{b'}+(\frac{4-6\lambda}{m}-\epsilon)G+(4\lambda-5+\epsilon)X_{b'}+\epsilon M'$ with $\epsilon$ a small rational number $0<\epsilon\ll 1$. 

 Let $\Delta=D_{|X_{b'}}$. Observe that $\j(X_b,F_{|X_b})=\j(X_b,6\Delta)$  is trivial for $\lambda < \frac{1}{864}$( and this follows as in the proof of Claim 1 and Corollary 2.35 of \cite{kol4}). We can choose such lambda for $\beta \ge 3266$. 
It follows that
\begin{eqnarray*}
H^0(X_b,\o_{X_b}(5K_{X_b})\otimes\j(X_b,6\Delta))&\cong& H^0(X_b,\o_{X_b}(5K_{X_b})).
\end{eqnarray*}
Consider the exact sequence 
\footnotesize$$
0\la \ox(5\k-X_b)\otimes\j(F)\la\ox(5\k)\otimes\j(F)\la\o_{X_b}(5K_{X_b})\otimes\j(F)_{|X_b}\la 0.
$$
\normalsize
Since $5\k-X_b-F\equiv \epsilon A'$ by Nadel vanishing (Theorem \ref{nadel}) we have a surjection
$$
H^0(X,\ox(5\k)\otimes\j(F))\la H^0(X_b,\o_{X_b}(5K_{X_b})\otimes\j(F)_{|X_b})\la 0.
$$
Observing that 
\begin{eqnarray*}
H^0(X_b,\o_{X_b}(5K_{X_b}))&\cong& H^0(X_b,\o_{X_b}(5K_{X_b})\otimes\j(X_b,F_{|X_b}))\\
&\subseteq&H^0(X_b,f_*\o_{X_b}(5K_{X_b})\otimes\j(F)_{|X_b}),
\end{eqnarray*}
and since we have sections in $\o_{X_b}(5K_{X_b})$ that separate $x$ and $y$ , we get  such a section in $\ox(5\k)$.
\qed

\bibliographystyle{plain}
\bibliography{bibl}

\end{document}